\documentclass[12pt,draft]{amsart}

\makeatletter

\theoremstyle{plain}
\newtheorem{thm}{Theorem}[section]
\numberwithin{equation}{section} 
\numberwithin{figure}{section} 
\theoremstyle{plain}
\newtheorem*{thm*}{Theorem}
\theoremstyle{plain}
\newtheorem{cor}[thm]{Corollary} 
\theoremstyle{plain}
\newtheorem*{cor*}{Corollary}

\theoremstyle{plain}
\newtheorem{lem}[thm]{Lemma} 
\theoremstyle{plain}
\newtheorem{prop}[thm]{Proposition} 
\theoremstyle{definition}
\newtheorem{defn}[thm]{Definition}
\theoremstyle{remark}
\newtheorem{rem}[thm]{Remark}
\theoremstyle{remark}

\theoremstyle{remark}

\theoremstyle{remark}
\newtheorem{question}[thm]{Question}
\theoremstyle{definition}

\theoremstyle{remark}
  \newtheorem*{acknowledgement*}{Acknowledgement}

\theoremstyle{plain}

\theoremstyle{plain}
\theoremstyle{plain}
\theoremstyle{plain}
\theoremstyle{definition}

\theoremstyle{remark}

\theoremstyle{remark}

\theoremstyle{remark}

\theoremstyle{plain}

\newcommand{\id}{\operatorname{id}}

\newcommand{\N}{\mathbb{N}}

\newcommand{\e}{\varepsilon}

\usepackage{amssymb} \usepackage{amscd}

\makeatother

\begin{document}

\title{Isomorphism of Hilbert modules over stably finite C$^*$-algebras}

\author{Nathanial P. Brown}
\address{Department of Mathematics, Penn State University,
State College, PA, 16802, USA}
\email{nbrown@math.psu.edu}
\author{Alin Ciuperca}
\address{Fields Institute, 222 College Street,  Toronto, Ontario, Canada, M5T 3J1}
\email{ciuperca@math.toronto.edu}

\thanks{N.B. was partially supported by DMS-0554870; A.C. was partially supported by Fields Institute.}

\begin{abstract}  It is shown that if $A$ is a stably finite C$^*$-algebra and $E$ is a countably generated Hilbert $A$-module, then $E$ gives rise to a compact element of the Cuntz semigroup if and only if $E$ is algebraically finitely generated and projective. It follows that if $E$ and $F$ are equivalent in the sense of Coward, Elliott and Ivanescu (CEI) and $E$ is algebraically finitely generated and projective, then $E$ and $F$ are isomorphic.  In contrast to this, we exhibit two CEI-equivalent Hilbert modules over a stably finite C$^*$-algebra that are not isomorphic.
\end{abstract}

\maketitle

\section{Introduction}

In \cite{cowelliottiv} a new equivalence relation -- we'll call it \emph{CEI equivalence} -- on Hilbert modules was introduced.  In general CEI equivalence is weaker than isomorphism, but it was shown that if $A$ has stable rank one, then it is the same as isomorphism (\cite[Theorem 3]{cowelliottiv}).  Quite naturally, the authors wondered whether their result could be extended to the stably finite case.  Unfortunately, it can't.  In Section \ref{sec:counterexample}, we give examples of Hilbert modules over a stably finite C$^*$-algebra which are CEI-equivalent, but not isomorphic. On the other hand, we show in Section \ref{sec:main} that CEI equivalence amounts to isomorphism when restricted to ``compact" elements of the Cuntz semigroup, in the stably finite case. 

\vspace{4mm}

\noindent\textbf{Acknowledgments:}  We thank George Elliott, Francesc Perera, Leonel Robert, Luis Santiago, Andrew Toms and Wilhelm Winter for valuable conversations on topics related to this work.

\section{Definitions and Preliminaries}

Throughout this note all C$^*$-algebras are assumed to be separable and all Hilbert modules are assumed to be right modules and countably generated.  We will follow standard terminology and notation in the theory of Hilbert modules (see, for example, \cite{lance}).  In particular, $\mathcal{K}$ denotes the compact operators on $\ell^2(\N)$, while $\mathcal{K}(E)$ will denote the ``compact" operators on a Hilbert module $E$.

For the reader's convenience, we recall a few definitions that are scattered throughout \cite{cowelliottiv}.

\begin{defn}
\label{defn:compactcontain}
If $E \subset F$ are Hilbert $A$-modules, we say $E$ is \emph{compactly contained in} $F$
if there exists a self-adjoint $T \in \mathcal{K}(F)$ such that $T|_E = \id_E$. In this situation we write $E \subset \subset F$.
\end{defn}

Note that $E \subset \subset E$ if and only if $\mathcal{K}(E)$ is unital; it can be shown that this is also equivalent to $E$ being algebraically finitely generated and projective (in the purely algebraic category of right $A$-modules) -- see the proof of \cite[Corollary 5]{cowelliottiv} (this part of the proof did not require the assumption of stable rank one.).

\begin{defn} We say a Hilbert $A$-module $E$ is \emph{CEI subequivalent} to another Hilbert $A$-module $F$ if every compactly contained submodule of $E$ is isomorphic to a compactly contained submodule of $F$.

We say $E$ and $F$ are \emph{CEI equivalent} if they are CEI subequivalent to each other -- i.e., a third Hilbert $A$-module $X$ is isomorphic to a compactly contained submodule of $E$ if and only if $X$ is isomorphic to a compactly contained submodule of $F$.
\end{defn}

\begin{defn}  We let $Cu(A)$ denote the set of Hilbert $A$-modules, modulo CEI equivalence.  The class of a module $E$ in $Cu(A)$ will be denoted $[E]$.
\end{defn}

It turns out that $Cu(A)$ is an abelian semigroup with $[E] + [F] := [E\oplus F]$. (Note: it isn't even obvious that this is well defined!) Moreover $Cu(A)$ is partially ordered -- $[E] \leq [F] \Longleftrightarrow$ $E$ is CEI subequivalent to $F$ -- and every increasing sequence has a supremum (i.e., least upper bound).  See \cite[Theorem 1]{cowelliottiv} for proofs of these facts.

\begin{defn} An element $x \in Cu(A)$ is \emph{compact} (in the order-theoretic sense) if for every increasing sequence $\{ x_n \} \subset Cu(A)$ with $x \leq \sup_n x_n$ there exists $n_0 \in \N$ such that $x \leq x_{n_0}$.
\end{defn}

For a unital C$^*$-algebra $A$, \emph{stable finiteness} means that for every $n \in \N$, $M_n(A)$ contains no infinite projections.  In the nonunital case there are competing definitions, but it seems most popular to say $A$ is stably finite if the unitization $\tilde{A}$ is stably finite, so this is the definition we will use.



\section{Main Results}
\label{sec:main}

The proof of our first lemma is essentially contained in the proof of \cite[Corollary 5]{cowelliottiv}.

\begin{lem}
\label{lem:equality}
Assume $E\subset \subset F$ is a compact inclusion of Hilbert $A$-modules.  If $E \cong F$ then either $E = F$ or $A\otimes \mathcal{K}$ contains a scaling element (in the sense of \cite{BC}).  If $A$ is stably finite, then $A\otimes \mathcal{K}$ cannot contain a scaling element; hence, in this case, $E \cong F$ if and only if $E = F$
\end{lem}

\begin{proof} Assume $E$ is properly contained in $F$; we'll show $A\otimes \mathcal{K}$ contains a scaling element. Let $v\colon F \to E$ be an isomorphism and $T \in \mathcal{K}(F)$ be a positive operator such that $T|_E = \id_E$. As observed in \cite{cowelliottiv}, the map $vT$ is adjointable -- i.e.\ defines an element of $\mathcal{L}(F)$ -- and, in fact, is compact. (This assertion is readily checked whenever $T$ is a ``finite-rank" operator). Moreover, a calculation shows that $(vT)^*|_E = Tv^{-1}$.  It is also worth noting that $T(vT) = vT$, since $T|_E = \id_E$ and $vT(F) \subset E$.

The scaling element we are after is $x = vT$.  Indeed, one checks that $x^* x = T^2$; hence,  $(x^* x)(xx^*) = T^2(vT)(vT)^* = (vT)(vT)^* = xx^*$.  Finally, we must see why $xx^* \neq x^* x$.  But if $xx^* = x^* x$, then $T^2 = (vT)(vT)^*$ and thus $T^2(F) \subset vT(F) \subset E$.  It follows that $T^2$ is a self-adjoint projection onto $E$ (since $T^2|_E = \id_E$, too), and hence $x = vT$ is a partial isometry whose support and range coincide with $E$. But this is impossible because $T = T^2$ (since $T \geq 0$), so $vT(F) \subsetneqq E$ (since $T(F) = E \subsetneqq F$).

We've shown that if $E \subsetneqq F$, then $\mathcal{K}(F)$ contains a scaling element.  But Kasparov's stabilization theorem provides us with an inclusion $\mathcal{K}(F) \subset A\otimes \mathcal{K}$, so the proof of the first part is complete.

In the case that $A$ is stably finite, it is well known to the experts that $A\otimes \mathcal{K}$ cannot contain a scaling element.  Indeed,  if it did, then \cite[Corollary 4.4]{BC} implies that $M_n(A)$ contains a scaling element, for some $n \in \N$. But it was shown in \cite{BC} that the unitization $\widetilde{M_n(A)}$ would then have an infinite projection. However, there is a natural embedding $\widetilde{M_n(A)} \subset M_n(\tilde{A})$, which contradicts the assumption of stable finiteness.
\end{proof}

Note that the canonical Hilbert module $\ell^2(A)$ is isomorphic to lots of (non-compactly contained) proper submodules.

\begin{prop}
\label{prop}
Let $E$ be a Hilbert $A$-module such that $[E]$ is compact in $Cu(A)$. Then either $E \subset \subset E$ or $A\otimes \mathcal{K}$ contains a scaling element.
\end{prop}

\begin{proof} Let $h \in \mathcal{K}(E)$ be strictly positive.  If $0$ is an isolated point in the spectrum $\sigma(h)$, then functional calculus provides a projection $p \in \mathcal{K}(E)$ such that $p = \id_E$;   so $E \subset \subset E$, in this case. If $0 \in \sigma(h)$ is not isolated, then, again using functional calculus, we can find $E_1 \subset \subset E_2 \subset \subset E_3 \cdots \subset \subset E$ such that $\cup_i E_i$ is dense in $E$ and $E_i \subsetneqq E_{i+1}$ for all $i \in \N$.

Since $[E]$ is compact, there exists $i$ such that $[E_i] = [E]$.  Since $E_{i+1} \subset \subset E$, $E_{i+1}$ is isomorphic to a compactly contained submodule of $E_i$ and this isomorphism restricted to $E_i$ maps onto a \emph{proper} submodule of $E_i$ (since $E_i \subsetneqq E_{i+1}$).  Thus $E_i$ is isomorphic to a proper compactly contained submodule of itself.  Hence, by Lemma \ref{lem:equality}, $A\otimes \mathcal{K}$ contains a scaling element.
\end{proof}

\begin{cor}  Let $A$ be stably finite and $E$ be a Hilbert $A$-module.  Then $[E] \in Cu(A)$ is compact if and only if $E \subset \subset E$. In particular, if $[E]$ is compact and $[E] \leq [F]$, then $E$ is isomorphic to a compactly contained submodule of $F$.
\end{cor}

\begin{proof}  The ``only if" direction is immediate from the previous proposition.  So assume $E \subset \subset E$ and let $[F_n] \in Cu(A)$ be an increasing sequence such that $[E] \leq [F] := \sup [F_n]$. By definition, $E$ is then isomorphic to a compactly contained submodule $E' \subset \subset F$. In the proof of \cite[Theorem 1]{cowelliottiv} it is shown that if $E' \subset \subset F$ and $[F] = \sup [F_n]$, then there is some $n \in \N$ such that $[E'] \leq [F_n]$.  Since $[E] = [E']$, the proof is complete.
\end{proof}

\begin{cor}
\label{cor:isom}
Let $A$ be stably finite and $E, F$ be Hilbert $A$-modules.  If $[E]= [F] \in Cu(A)$ is compact, then $E \cong F$.  In particular, if $[E]= [F]$ and $E$ is algebraically finitely generated and projective, then $[E] \in Cu(A)$ is compact; hence, $E \cong F$.
\end{cor}

\begin{proof}  Assume $[E] = [F]$ is compact.  Then $E \subset \subset E$ and $F \subset \subset F$, by the previous corollary.  Hence there exist isomorphisms $v\colon F \to F' \subset \subset E$ and $u\colon E \to E' \subset \subset F$.  It follows that $F \cong u(v(F)) \subset \subset F$, which, by Lemma \ref{lem:equality}, implies that $u(v(F)) = F$.  Hence $u$ is surjective, as desired.

As mentioned after Definition \ref{defn:compactcontain},  if $E$ is algebraically finitely generated and projective, then $E \subset \subset E$, which implies $[E]$ is compact (as we've seen).  
\end{proof}

In the appendix of \cite{cowelliottiv} it is shown that $Cu(A)$ is isomorphic to the classical Cuntz semigroup $W(A\otimes \mathcal{K})$.  When $A$ is stable, the isomorphism $W(A) \to Cu(A)$ is very easy to describe: the Cuntz class of $a \in A_+$ is sent to $H_a := \overline{aA}$ (with its canonical Hilbert $A$-module structure).

\begin{thm}
\label{thm:main}
Let $A$ be a stable, finite C$^*$-algebra, $a \in A_+$ and $H_a = \overline{aA}$. The following are equivalent:
\begin{enumerate}
\item $H_a$ is algebraically finitely generated and projective;

\item $[H_a] \in Cu(A)$ is compact;

\item $\sigma(a) \subset \{0\} \cup [\e, \infty)$ for some $\e > 0$;

\item $\langle a \rangle = \langle p \rangle \in W(A)$ for some projection $p\in A$.
\end{enumerate}
\end{thm}

\begin{proof}  The implication $(1) \Longrightarrow (2)$ was explained above.

$(2) \Longrightarrow (3)$: Let $a_\e = (a-\e)_+$.  Then $H_{a_\e} \subset \subset H_a$ and $\cup_\e H_{a_\e}$ is dense in $H_a$. Since $[H_a] \in Cu(A)$ is compact, there exists $\e > 0$ such that $[H_a] = [H_{a_\e}]$.  Corollary \ref{cor:isom} implies that $H_a \cong H_{a_\e}$; thus $H_a = H_{a_\e}$, by Lemma \ref{lem:equality}.  It follows that $\sigma(a) \subset \{0\} \cup [\e, \infty)$, because otherwise functional calculus would provide a nonzero element $b \in C^*(a)$ such that $0 \leq b \leq a$ (so $b \in H_a$) and $a_\e b = 0$ (so $b \notin H_{a_\e}$), which would contradict the equality $H_a = H_{a_\e}$.

$(3) \Longrightarrow (4)$ is a routine functional calculus exercise.

$(4) \Longrightarrow (1)$: Assume $\langle a \rangle = \langle p \rangle \in W(A)$.  Since $pA$ is singly generated and algebraically projective, Corollary \ref{cor:isom} implies $H_a$ is isomorphic to $pA$.
\end{proof}

The equivalence of $(3)$ and $(4)$ above generalizes Proposition 2.8 in \cite{PT}.

\begin{cor} If $A$ is stably finite, then $A\otimes \mathcal{K}$ has no nonzero projections if and only if $Cu(A)$ contains no compact element.
\end{cor}

\section{A Counterexample}
\label{sec:counterexample}

Now let us show that  if $A$ is stably finite and $E,F$ are Hilbert $A$-modules such that $[E] = [F]$, then it need not be true that $E$ and $F$ are isomorphic.  Let $A = C_0(0,1] \otimes \mathcal{O}_3 \otimes \mathcal{K}$, where $\mathcal{O}_3$ is the Cuntz algebra with three generators. Voiculescu's homotopy invariance theorem (cf.\ \cite{dvv}) implies that $A$ is quasidiagonal,  hence stably finite. Let $p, q \in \mathcal{O}_3 \otimes \mathcal{K}$ be two nonzero projections which are \emph{not} Murray-von Neumann equivalent. If $x \in C_0(0,1]$ denotes the function $t \mapsto t$, then we define $f_p = x \otimes p$ and $f_q = x \otimes q$ in $A$. Since $A$ is purely infinite in the sense of \cite{KR} and the ideals generated by $f_p$ and $f_q$ coincide, it follows that $[\overline{f_p A}] = [\overline{f_q A}] \in Cu(A)$. We claim that the modules $\overline{f_p A}$ and $\overline{f_q A}$ are not isomorphic.

Indeed, if they were isomorphic, then we could find $v \in A$ such that $v^* v = f_p$ and $\overline{vv^*A} = \overline{f_q A}$. (See \cite[Lemma 3.4.2]{ciuperca}; if $T\colon \overline{f_p A}\to\overline{f_q A}$ is an isomorphism, then $v = T(f_p^{1/2})$ has the asserted properties.) Letting $\pi\colon A \to \mathcal{O}_3 \otimes \mathcal{K}$ be the quotient map corresponding to evaluation at $1 \in (0,1]$, it follows that $\pi(v)^* \pi(v) = p$ and $\overline{\pi(v) \pi(v)^* (\mathcal{O}_3 \otimes \mathcal{K})} = \overline{q (\mathcal{O}_3 \otimes \mathcal{K})}$.  Since $\pi(v) \pi(v)^*$ is a projection whose associated hereditary subalgebra agrees with the hereditary subalgebra generated by $q$, it follows that $\pi(v) \pi(v)^* = q$ (since both projections are units for the same algebra). This contradicts the assumption that $p$ and $q$ are not Murray-von Neumann equivalent, so $\overline{f_p A}$ and $\overline{f_q A}$ cannot be isomorphic.

\section{Questions and Related Results}

If the following question has an affirmative answer, then the proof of \cite[Corollary 5]{cowelliottiv} would show that $A$ has real rank zero if and only if the compacts are ``dense" in $Cu(A)$.

\begin{question} Can Corollary \ref{cor:isom} be extended to the ``closure" of the compact elements?  That is, if $A$ is stably finite and $E$ and $F$ are Hilbert A-modules such that $[E]=[F] = \sup [C_n]$ for an increasing sequence of compact elements $[C_n]$, does it follow that $E\cong F$?
\end{question}

The next question was raised in \cite{cowelliottiv}, but we repeat it because the modules in Section \ref{sec:counterexample} are not counterexamples -- they mutually embed into each other.   (To prove this, use the fact that $p$ is Murray-von Neumann equivalent to a subprojection of $q$, and vice versa.)

\begin{question}  Are there two Hilbert modules $E$ and $F$ such that $[E] = [F]$, but $F$ is not isomorphic to a submodule of $E$?
\end{question}

\begin{question}  If $x \in Cu(A)$ is compact, is there a projection $p \in A\otimes \mathcal{K}$ such that $x = \langle p \rangle$?
\end{question}

Of course, in the stably finite case the results of Section \ref{sec:main} tell us that much more is true, but for general C$^*$-algebras we don't know the answer to this question.  However, we can give an affirmative answer in some interesting cases, as demonstrated below.   First, a definition.

\begin{defn} An element $x \in Cu(A)$ will be called \emph{infinite} if $x+y=x$ for some non-zero $y\in Cu(A)$. Otherwise, $x$ will be called \emph{finite}.
\end{defn}

Note that $[\ell^2(A)] \in Cu(A)$ is always infinite.

\begin{lem}
\label{lem:unique}
If $A$ is simple, then $[\ell^2(A)] \in Cu(A)$ is the unique infinite element.
\end{lem}

\begin{proof} Assume $[E] + [F] = [E]$ for some nonzero Hilbert $A$-module $F$.  Adding $[F]$ to both sides, we see that $[E] + 2[F] = [E]$; repeating this, we have that $[E] + k[F] = [E]$ for all $k \in \N$. By uniqueness of suprema, it follows that $[E] + [\ell^2(F)] = [E]$ (cf.\ \cite[Theorem 1]{cowelliottiv}).  Since $A$ is simple, $F$ is necessarily full and hence $\ell^2(F) \cong \ell^2(A)$ (\cite[Proposition 7.4]{lance}).  Thus $$[E] = [E] + [\ell^2(F)] = [E \oplus \ell^2(A)] = [\ell^2(A)],$$ by Kasparov's stabilization theorem.
\end{proof}

In the proof of the following lemma, we use the operator inequality $$xbx^*+ y^*by\geq xby + y^*bx^*,$$ for any $b$ in $A^+$, and $x, y\in A$. (Which follows from the fact that $(x-y^*)b(x-y^*)^*\geq0$.)

\begin{lem}\label{algsimple} Let $A$ be a stable algebraically simple C*-algebra.
\begin{enumerate}
\item For any non-zero $x\in Cu(A)$ there exists $n\in \N$ such that $nx=[A]$.

\item There exists a projection $q\in A$ such that $[A]=[qA]$. In particular, $[A]$ is a compact element of the Cuntz semigroup $Cu(A)$.
\end{enumerate}
\end{lem}

\begin{proof} It will be convenient to work in the original positive-element picture of the Cuntz semigroup.  Our notation is by now standard (cf.\ \cite{PT}).

Proof of (1): Let $x=[\overline{bA}]$ for some $0\neq b\in A^+$ and let $a\in A$ be a strictly positive element.
(Stability implies that every right Hilbert $A$-module is isomorphic to a closed right ideal of $A$.) Since $A$
is algebraically simple, one can find $x_1,\ldots, x_n, y_1, \ldots, y_n \in A$ such that $a=\sum_{i=1}^k
x_iby_i$. Thus,
\begin{align*}
 a\sim 2a=a+a^* &=\sum_{i=1}^k (x_iby_i +  y_i^*bx_i^*)\\
 & \leq \sum_{i=1}^k (x_ibx_i^*+y_i^*by_i)\\
& \lesssim x_1bx_1^*\oplus y_1^*by_1 \oplus \cdots x_kbx_k^*\oplus y_k^*by_k\\
& \lesssim b\oplus b\oplus \cdots \oplus b,
\end{align*}
where the last sum has $n=2k$ summands.

Since $A$ is stable, one can embed the Cuntz algebra $O_n$ in the multiplier algebra $M(A)$. This gives us isometries $s_1,\cdots, s_n\in M(A)$ with orthogonal ranges. Set $b_i'=s_ibs_i^*$ and note that $b_i'\sim b$ and $b_i'\perp b_j'$. Moreover, $a\lesssim b_1'+\cdots +b_n'\lesssim a$ (since $a$ is  strictly positive, it Cuntz-dominates any element of $A$). Therefore, $\langle a \rangle = n\langle b\rangle = nx$, or equivalently, $[A] = nx$.

Proof of (2): Since $A$ is stable and algebraically simple,  \cite[Theorem 3.1]{BC} implies $A$ has a non-zero projection $p$. As above, we can find orthogonal projections $p_1,\ldots, p_n \in A$ such that $p_i\sim p$ and $\langle p_1+\cdots +p_n \rangle = n\langle p \rangle = [A]$. Defining $q = p_1+\cdots +p_n$, we are done.
\end{proof}

We'll also need a consequence of the work in Section \ref{sec:main}.

\begin{prop}
\label{stablyfinite}
If $A$ is stable, $\langle a \rangle \in W(A) = Cu(A)$ is compact and $0 \in \sigma(a)$ is not an isolated point, then $A$ contains a scaling element and $\langle a \rangle$ is infinite.
\end{prop}

\begin{proof}  Assume $A$ contains no scaling element.  Since $\langle a \rangle$ is compact, Proposition \ref{prop} implies that $H_a \subset \subset H_a$.  As in the proof of $(2) \Longrightarrow (3)$ in Theorem \ref{thm:main}, there exists $\e > 0$ such that $[H_a] = [H_{a_\e}]$ and hence $H_a$ is isomorphic to a compactly contained submodule $E$ of $H_{a_\e}$. Lemma \ref{lem:equality} implies $E = H_a$, so $H_{a_\e} = H_a$ too. As we've seen, this implies $\sigma(a) \subset \{0\} \cup [\e, \infty)$, contradicting our hypothesis; hence, $A$ contains a scaling element.

To prove the second assertion, choose $\e > 0$ such that $[H_a] = [H_{a_\e}]$. Since $0 \in \sigma(a)$ is not isolated, we can find a nonzero positive function $f \in C_0(0,\|a\|]$ such that $f(t) = 0$ for all $t \geq \e$.  Thus $f(a) + (a-\e)_+ \precsim a$ and $f(a) (a-\e)_+ = 0$. It follows that $$[H_{f(a)}] + [H_a] = [H_{f(a)}] + [H_{a_\e}] \leq [H_a]$$ and thus $[H_a]$ is infinite.
\end{proof}

\begin{thm}  Let $x \in Cu(A)$ be compact.
\begin{enumerate}
\item If $A$ is simple, then there exists a projection $p \in A\otimes \mathcal{K}$ such that $x = \langle p \rangle$.

\item If $x$ is finite, then there exists a projection $p \in A\otimes \mathcal{K}$ such that $x = \langle p \rangle$.
\end{enumerate}
\end{thm}

\begin{proof}  In both cases we may assume $A$ is stable.

Proof of (1):  Fix a nonzero positive element $a\in A$ such that $x = [H_a]$.  If $0 \in \sigma(a)$ is an isolated point, then functional calculus provides us with a Cuntz equivalent projection, and we're done. Otherwise Proposition \ref{stablyfinite} tells us that $x$ is infinite and $A$ contains a scaling element. By simplicity and Lemma \ref{lem:unique}, we have that $x = [\ell^2(A)] = [A]$ (by stability).  Moreover, the existence of a scaling element ensures that $A$ is algebraically simple (see \cite[Theorem 1.2]{BC}).  Hence part (2) of Lemma \ref{algsimple} provides the desired projection.

Proof of (2):  Choose $a \in A_+$ such that $x = \langle a \rangle$.  Since $x$ is finite, Proposition \ref{stablyfinite} implies $0 \in \sigma(a)$ is an isolated point, so we're done.
\end{proof}

\begin{rem}  It is possible to improve part (2) of the theorem above.  Namely, it is shown in \cite{ciuperca} that if $x \in Cu(A)$ is compact and there is no \emph{compact} element $y \in Cu(A)$ such that $x = x + y$, then there exists a projection $p \in A\otimes \mathcal{K}$ such that $x = \langle p \rangle$.
\end{rem}

\end{document}